\documentclass{article}

\usepackage{mystyle}
\usepackage{wasysym}
\usepackage{mathtools}
\begin{document}
\title{A note on cardinal preserving embeddings}
\author{Gabriel Goldberg\\ Evans Hall\\ University Drive \\ Berkeley, CA 94720}
\maketitle
\section{Introduction}
An elementary embedding between two transitive models of set theory is \textit{cardinal preserving}
if its domain and target models have the same class of cardinals.
The purpose of this note is to extend some ideas due to Caicedo and Woodin, leading
to a proof that the critical point \(\kappa\) of a cardinal preserving embedding
from the universe into an inner model 
is strongly compact in \(V_\gamma\) for some inaccessible cardinal \(\gamma > \kappa\).

In particular, we have the following consistency result:
\begin{repthm}{ConThm}
The existence of a cardinal preserving embedding from the universe into an inner model
implies the consistency of \textnormal{ZFC + }there is a strongly compact cardinal. 
\end{repthm}
We also show that the existence of a cardinal preserving
embedding from the universe into \(M\) implies much stronger resemblance properties of
\(V\) and \(M\) than are immediately apparent, especially assuming the Singular Cardinals Hypothesis:
\begin{repthm}{SCHThm}[SCH]
    Suppose there is a cardinal preserving elementary embedding
    from the universe into an inner model \(M\). Then for any ordinal \(\alpha\), \(\cf^M(\alpha) = \cf(\alpha)\)
    and for any cardinal \(\rho\), \((2^\rho)^M = 2^\rho\).
\end{repthm}
\section{Some proofs}
If \(j : V\to M\) is an elementary embedding, then the critical sequence \(\kappa_n(j)\) is defined by recursion,
setting \(\kappa_0(j) = \crit(j)\) and \(\kappa_{n+1}(j) = j(\kappa_n(j))\). Finally,
\(\kappa_\omega(j) = \sup_{n < \omega}\kappa_n(j)\) denotes the supremum of the critical sequence of \(j\).

The following is straightforward using an observation of Caicedo:
\begin{prp}\label{CompactPrp}
    Suppose there is a cardinal preserving embedding \(j : V\to M\).
    Let \(\kappa = \crit(j)\) and let \(\lambda = \kappa_\omega(j)\).
    Then \(\kappa\) is \(\lambda\)-strongly compact.
    \begin{proof}
        By Ketonen's Theorem \cite{Ketonen}, it suffices to show that every regular cardinal
        in the interval \([\kappa,\lambda]\) carries a \(\kappa\)-complete uniform ultrafilter.
        By Pr\'ikry's Theorem \cite{Prikry}, it suffices to show that every successor
        cardinal in the interval \([\kappa,\lambda]\) carries a \(\kappa\)-complete uniform ultrafilter.

        Suppose \(\delta\) is an ordinal with \(\kappa\leq \delta < \lambda\). We claim there is a \(\kappa\)-complete
        uniform ultrafilter on \(\delta^+\). First, since \(\kappa\leq \delta^+ < \lambda\), \(\delta^+ < j(\delta^+)\).
        Also, note that \(j(\delta^+) = j(\delta)^+\) since \(M\)
        correctly computes successor cardinals. 
        This means that the cardinal \(j(\delta^+)\) is regular, and in particular it has cofinality
        strictly larger than \(\delta^+\). As a consequence, \(\sup j[\delta^+] < j(\delta^+)\).
        Therefore one can derive an ultrafilter on \(\delta^+\) from \(j\) using \(\sup j[\delta^+]\),
        and this ultrafilter is clearly uniform; it is \(\kappa\)-complete since \(\crit(j) = \kappa\).
    \end{proof}
\end{prp}
The proposition does not obviously yield the consistency of a strongly compact cardinal from
a cardinal preserving embedding since it is not clear that
the interval \((\kappa,\lambda)\) contains an inaccessible cardinal. In fact, it is not
even clear that \(2^\kappa < \lambda\).
We will prove that \(\lambda\) is a limit of inaccessible cardinals,
and in fact \(\kappa_n(j)\) is inaccessible for every \(n < \omega\).

The key lemma is not hard, but it seems interesting in its own right.
The main concept involved is the \textit{tightness function} of an elementary embedding.
\begin{defn}
    If \(j : V\to M\) is an elementary embedding,
    then for any cardinal \(\lambda\), \(\tr{j}{\lambda}\) denotes the
    least \(M\)-cardinality of a set in \(M\) that covers \(j[\lambda]\).
\end{defn}
Note that for any set \(X\), \(\tr{j}{|X|}\) is the least \(M\)-cardinality of a set in \(M\) that covers
\(j[X]\). The following theorem, implicit in the work of Ketonen \cite{Ketonen}, 
appears explicitly in the author's thesis \cite[Theorem 7.2.12]{UA}.
\begin{thm}[Ketonen, \cite{Ketonen}]\label{KetonenTightness}
    If \(j : V\to M\) is an elementary embedding and \(\delta\) is a regular cardinal,
    then \(t_j(\delta) = \cf^M(\sup j[\delta])\).\qed
\end{thm}
\begin{lma}\label{KeyLma}
    Suppose \(j : V\to M\) is an elementary embedding with critical point \(\kappa\)
    and \(\lambda\) is a singular cardinal of cofinality \(\iota < \kappa\).
    Then for some \(M\)-cardinal \(\eta\) of \(M\)-cofinality \(\iota\),
    \(\eta < \tr{j}{\lambda^\iota} \leq (\eta^{\iota})^M\).
    \begin{proof}
        Let \(U\) be the ultrafilter on \(\lambda^+\) derived from \(j\)
        using \(\sup j[\lambda^+]\). Let \(i : V\to N\) be the ultrapower associated to \(U\). 

        Let \(\langle \delta_\alpha : \alpha < \iota\rangle\) be an increasing sequence of
        regular cardinals cofinal in \(\lambda\).
        Let \(\delta_\alpha' = \tr{i}{\delta_\alpha}\). Then since \(N\) is closed under \(\iota\)-sequences,
        \(\tr{i}{\lambda} = \sup \delta_{\alpha}'\). Moreover, \cref{KetonenTightness} easily implies
        that \(\delta_\alpha' < \delta_\beta'\) whenever \(\alpha < \beta\).
        Therefore \(t_i(\lambda)\) has cofinality \(\iota\) (in \(N\)).
        Let \(\nu = t_i(\lambda)\).
        Fix \(A\in N\) covering \(i[\lambda]\) with \(|A|^N = \nu\).
        Since \(N\) is closed under \(\iota\)-sequences and \(\crit(j) > \iota\),
        \(i[{}^\iota\lambda]\subseteq ({}^\iota A )^N\).
        It follows that \(\tr{i}{\lambda^\iota}\leq \nu^\iota\).

        Let \(k : N\to M\) be the factor embedding, and set \(\eta = k(\nu)\).
        Trivially, \(\tr{j}{\lambda^\iota} \leq (\eta^{\iota})^M\).
        Since \(k(\sup i[\lambda^+]) = \sup j[\lambda^+]\),
        \cref{KetonenTightness} implies that \(k(\tr{i}{\lambda^+}) = \tr{j}{\lambda^+}\).
        Therefore since \(\nu = \tr{i}{\lambda} < \tr{i}{\lambda^+}\),
        \(\eta < \tr{j}{\lambda^+} \leq \tr{j}{\lambda^\iota}\). This finishes the proof.
    \end{proof}
\end{lma}
Here is one intriguing consequence, which is very similar to a consequence of the Ultrapower Axiom 
\cite[Lemma 8.4.14]{UA}
originally proved in the analysis of cardinal preservation under UA.
\begin{thm}[SCH]\label{SCHLma}
    Suppose \(j : V\to M\) is an elementary embedding and 
    \(\delta\) is the successor of a singular cardinal of cofinality \(\iota < \crit(j)\).
    Then in \(M\), the cofinality of \(\sup j[\delta]\) is the 
    successor of a singular cardinal of cofinality \(\iota\).\qed
\end{thm}
We will avoid the SCH in our main theorem by citing a local form of Solovay's Theorem.
\begin{thm}[{Solovay, \cite{Solovay}}]\label{SolovaySCH}
    If \(\kappa\) is \(\lambda\)-strongly compact, then for any singular
    cardinal \(\eta\) with \(\kappa\leq \eta \leq \lambda\),
    \(\eta^{\cf(\eta)} = \eta^+\cdot 2^{\cf(\eta)}\).\qed
\end{thm}
The author tentatively conjectures that the existence of a cardinal preserving embedding
of the universe of sets is inconsistent with ZFC. The existence of the partially cardinal
preserving embeddings involved in the hypotheses of the next few theorems, however, 
follows from the axiom \(I_2\), and therefore is very likely to be consistent.
\begin{thm}\label{CofPreserve}
    Suppose \(j : V\to M\) is an elementary embedding.
    Let \(\kappa = \crit(j)\), and let \(\lambda = \kappa_\omega(j)\).
    Assume \(j\restriction V_\lambda\) is cardinal preserving,
    or in other words \(\Card^M\cap \lambda = \Card^V \cap \lambda\).
    Then for any ordinal \(\alpha \leq \lambda\), \(\cf^M(\alpha) = \cf(\alpha)\).
    \begin{proof}
        First, fix a singular cardinal \(\rho\) of cofinality \(\omega\)
        with \(\rho \leq \lambda\).
        We will show that \(\tr{j}{\rho^+} = \rho^+\). 
        If \(j(\rho^{+}) = \rho^+\), this is trivial, so we may assume \(\kappa \leq \rho < \lambda\).
        By \cref{KeyLma}, there is a cardinal \(\eta\)
        of \(M\)-cofinality \(\omega\) such that 
        \(\eta < \tr{j}{\lambda^+}\leq (\eta^\omega)^M\).

        We claim \((\eta^\omega)^M = \eta^+\).
        Note that \(\kappa \leq \eta < \lambda\), so
        since \(\kappa\) is \(\lambda\)-strongly compact by \cref{CompactPrp},
        \(\eta^\omega = \eta^+\) by Solovay's Theorem.
        Therefore \((\eta^\omega)^M < \eta^{++}\), and so by cardinal correctness,
        \((\eta^\omega)^M \leq \eta^+\). By K\"onig's Theorem, \((\eta^\omega)^M \geq (\eta^+)^M = \eta^+\).
        It follows that \((\eta^\omega)^M = \eta^+\), as claimed.

        By \cref{KeyLma}, we can conclude that \(\eta < \tr{j}{\lambda^+}\leq \eta^+\),
        or in other words, \(\tr{j}{\rho^+} = \eta^+\).
        Hence \(\cf^{M}(\sup j[\rho^+]) = \eta^+\), and
        so \(\eta^+\) and \(\rho^+\) have the same cofinality. Since successor cardinals
        are regular, \(\eta^+ = \rho^+\).

        Now let \(\gamma < \lambda\) be a regular cardinal. We have 
        \[\gamma \leq \tr {j}{\gamma} < \tr{j}{\gamma^{+\omega+1}} = \gamma^{+\omega+1}\]
        where the final equality follows from the previous paragraphs
        with \(\rho = \gamma^{+\omega}\).
        Thus \(\tr{j}{\gamma} = \gamma^{+n}\) for some \(n < \omega\). But
        by \cref{KetonenTightness}, \(\tr{j}{\gamma}\) has cofinality \(\gamma\),
        and therefore \(\tr{j}{\gamma} =\gamma\).

        Next, suppose \(\delta < \lambda\) is a cardinal that is regular in \(M\).
        Assume towards a contradiction that \(\delta\) is singular in \(V\).
        By our work so far, \(\delta \leq \tr{j}{\delta} < \tr{j}{\delta^+} = \delta^+\),
        so \(\tr{j}{\delta} = \delta\), and hence \(\cf^M(\sup j[\delta]) = \delta\).
        But now let \(\iota = \cf(\delta)\).
        Let \(U\) be the ultrafilter on \(\iota\) derived from \(j\) using \(\sup j[\iota]\).
        Then \(\tr{j_U}{\iota} = \iota\).
        As a consequence, \(\cf^{M_U}(\sup j_U[\delta]) = \iota\).
        Now let \(k : M_U\to M\) be the factor embedding.
        Then \(k(\iota) = \iota\). Hence \(k\) is continuous
        at ordinals of cofinality \(\iota\),
        and in particular, \(k(\sup j_U[\delta]) = \sup j[\delta]\).
        Thus \(\cf^M(\sup j[\delta]) = k(\cf^{M_U}(\sup j_U[\delta])) = k(\iota) = \iota\),
        and this contradicts that \(\cf^M(\sup j[\delta]) = \delta\).

        Finally, we conclude the theorem. Fix an ordinal \(\alpha\).
        Then \(\cf^M(\alpha)\) is an \(M\)-regular cardinal \(\delta\),
        and so \(\cf(\alpha) = \cf(\delta) = \delta\).
    \end{proof}
\end{thm}
We now argue that cardinal preserving embeddings preserve the continuum function.
This extends an argument due to Woodin-Caicedo \cite{Caicedo}. 
\begin{lma}\label{TightContinuum}
    Suppose \(j :V \to M\) is an elementary embedding. If 
    \(\rho\) is a cardinal and \(\theta = \tr{j}{\rho}\)
    then \((2^\theta)^M \geq 2^\rho\).
    \begin{proof}
        Let \(A\in M\) be a cover of \(j[\rho]\) with \(|A|^M = \theta\).
        Note that \(|P^M(A)|^M = (2^\theta)^M\).
        But there is an injection \(f : P(\rho)\to P^M(A)\)
        defined by \(f(B) = j(B)\cap A\).
        So \(2^\rho = |P(\rho)| \leq |P^M(A)| \leq |P^M(A)|^M = (2^\theta)^M\).
    \end{proof}
\end{lma}

\begin{thm}\label{ContinuumPreserving}
    Suppose \(j : V\to M\) is an elementary embedding.
    Let \(\kappa = \crit(j)\), and let \(\lambda = \kappa_\omega(j)\).
    Assume \(j\restriction V_\lambda\) is cardinal preserving,
    or in other words \(\Card^M\cap \lambda = \Card^V \cap \lambda\).
    Then for any cardinal \(\rho < \lambda\), \((2^\rho)^M = 2^\rho\).
    \begin{proof}
        Since \(|(2^\rho)^M| = |P^M(\rho)|\), we have
        \((2^\rho)^M < (2^\rho)^+\). But \((2^\rho)^M\) is a cardinal,
        so it follows that \((2^\rho)^M \leq 2^\rho\).
        On the other hand, by \cref{CofPreserve},
        \(\tr{j}{\rho} = \rho\), and so by \cref{TightContinuum}, 
        \((2^\rho)^M\geq 2^\rho\).
    \end{proof}
\end{thm}

\begin{thm}\label{ConThm}
    The existence of a cardinal preserving embedding from the universe into an inner model
    implies the consistency of \textnormal{ZFC + }there is a strongly compact cardinal. 
    \begin{proof}
        If \(j : V\to M\) is cardinal preserving, then by a simple induction using \cref{ContinuumPreserving},
        \(\kappa_n(j)\) is strongly inaccessible for all \(n < \omega\).
        By \cref{CompactPrp}, for any \(n \geq 1\), 
        \(V_{\kappa_n(j)}\) satisfies ZFC + there is a strongly compact cardinal.
    \end{proof}
\end{thm}
\cref{ContinuumPreserving} has some surprising reflection consequences:
for example, the SCH must hold for all sufficiently large cardinals below the critical point \(\kappa\)
of a cardinal preserving embedding \(j : V\to M\). 
(It is open whether the SCH can fail everywhere below a strongly compact cardinal.)
Moreover, if GCH holds below \(\kappa\), it holds up to \(\kappa_\omega(j)\). 
Can \(\kappa\) be the least measurable cardinal?

The same arguments almost prove that any cardinal preserving embedding \(j: V\to M\) is cofinality preserving
and continuum function preserving,
but now there is a problem applying Solovay's Theorem that SCH holds above a strongly compact cardinal.
For example, letting \(\lambda = \kappa_\omega(j)\), it is not clear how large \(\lambda^\omega\) can be.
We do obtain: 
\begin{thm}[SCH]\label{SCHThm}
    Suppose there is a cardinal preserving elementary embedding
    from the universe into the inner model \(M\). Then for any ordinal \(\alpha\), \(\cf^M(\alpha) = \cf(\alpha)\)
    and for any cardinal \(\rho\), \((2^\rho)^M = 2^\rho\).\qed
\end{thm}
\bibliographystyle{plain}
\bibliography{Bibliography.bib}
\end{document}